\def\versiondate{25 Sep. 2011}
\input math.macros
\input Ref.macros

\checkdefinedreferencetrue
\continuousnumberingtrue
\continuousfigurenumberingtrue
\theoremcountingtrue
\sectionnumberstrue
\forwardreferencetrue
\citationgenerationtrue
\nobracketcittrue
\hyperstrue
\initialeqmacro

\def\gp{\Gamma}
\def\gpe{\gamma}
\def\gh{G}
\def\pc{p_{\rm c}}

\def\rr{{\scr R}}

\input\jobname.key
\bibsty{../../texstuff/myapalike}

\ifproofmode \relax \else\head{To appear in {\it Erg. Th. Dyn. Sys.}}
{Version of \versiondate}\fi 
\vglue20pt

\title{Fixed Price of Groups and Percolation}

\author{Russell Lyons}

\abstract{We prove that for every finitely generated
group $\gp$, at least one of the following holds:
(1) $\gp$ has fixed price; (2) each of its Cayley graphs $\gh$ has
infinitely many infinite clusters for some Bernoulli percolation on $\gh$.}

\bottomII{Primary 
22F10,	
60K35. 
}
{Measurable group actions, cost.}
{Research partially supported by NSF grant DMS-1007244.}

\bsection{Introduction}{s.intro}

Let $\gp$ be an infinite finitely generated group.
Consider an essentially free measure-preserving action of $\gp$ on a
standard probability space $(X, \mu)$. A \dfn{graphing} of this action is a
graph $(X, E)$, where $E \subseteq X \times X$ is measurable and symmetric
and such that for every $x \in X$, the vertices of the connected component
of $x$ are the same as the orbit $\gp x$.
The \dfn{cost} of such a graphing equals $(1/2) \int_X |\{\gpe \in \gp \st
(x, \gpe x) \in E\}| \,d\mu(x)$.
The \dfn{cost} of the action
is the infimum of the costs of its graphings.
See \ref b.Gaboriau:whatis/ for background information on cost.
Following \ref b.Gaboriau:cost/, we say that $\gp$ has \dfn{fixed price}
if all its (essentially free measure-preserving) actions have the same
cost. One of the major open questions of the
theory of cost is whether every group has fixed price. 
It is known that all amenable groups have fixed price with cost 1. Some
other restricted classes and examples of groups of fixed price are also
known; see Section VII of \ref b.Gaboriau:cost/.

Now let $\gh$ be a Cayley graph of $\gp$.
Bernoulli($p$) site \dfn{percolation} on $\gh$ is the subgraph of $\gh$
induced by a random subset of its vertices, where each vertex is in the
subset with probability $p$ independently of other vertices.
The connected components of this graph 
are called \dfn{clusters}.
If $\gp$ is amenable, then \ref b.BK:uni/ and \ref b.GKN:uni/ proved that
there is at most one infinite cluster a.s.
A major open conjecture of \ref b.pyond/ is the converse, that if $\gp$
is non-amenable, then for some interval of $p$, there are infinitely many
infinite clusters a.s.
By work of \ref b.BLPS:gip/, it is known that this is equivalent to the
existence of a single value of $p$ at which there are infinite many
infinite clusters a.s.
\ref b.Lyons:review/ observed that {\it every} Cayley graph has this
property when $\gp$ has cost $> 1$; no other groups are known such that all
its Cayley graphs have this property.
\ref b.PakSN:uniq/ proved that $\gp$ has {\it some} Cayley graph with this
property for {\it bond} percolation when $\gp$ is non-amenable, where in
bond percolation, it is the edges, not the vertices, that are kept with
probability $p$ independently.

Here we note that if $\gp$ does not have fixed price with cost 1, then for
each of its Cayley graphs $\gh$, there exists a $p$ where Bernoulli($p$)
percolation has infinitely many infinite clusters a.s.
We prove this for site percolation, but essentially
the same proof applies for bond percolation.
Thus, our result here includes that of \ref b.Lyons:review/, but actually
the method of proof (not written anywhere) is the same. What is new is the
combination of that method with a result of \ref b.AW:Bernoulli/.

\bsection{Proof}{s.proof}

Consider a measurable equivalence relation $\rr$ with countable equivalence
classes on a
standard probability space $(X, \mu)$. A \dfn{graphing} of $\rr$ is a
graph $(X, E)$, where $E \subseteq \rr$ is measurable and symmetric
and such that for every $x \in X$, the vertices of the connected component
of $x$ are the same as the equivalence class of $x$.
The \dfn{cost} of such a graphing equals $(1/2) \int_X |\{y \st
(x, y) \in E\}| \,d\mu(x)$.
The \dfn{cost} of $\rr$
is the infimum of the costs of its graphings.

\ref b.AW:Bernoulli/ proved that every Bernoulli action of $\gp$ has the
maximum cost among all essentially free measure-preserving actions of $\gp$.
Let $\gh$ be a Cayley graph of a non-amenable group $\gp$ with respect to a
generating set $S$.
Let $\theta(p)$ denote the probability that a given vertex of $\gh$ belongs
to an infinite cluster in Bernoulli($p$) percolation; this is the same for
all vertices.
Define $\pc(\gh) := \inf \big\{ p \in [0, 1] \st \theta(p) > 0\big\}$.
If there are no $p$ with infinitely many infinite clusters a.s., then for
all $p > \pc(\gh)$, there is a unique infinite cluster a.s.\ by a theorem
of \ref b.NS/. Furthermore, \ref
b.BLPS:gip/ proved that $\theta\big(\pc(\gh)\big) = 0$, whence by \ref
b.BK:conti/, $\lim_{p \downarrow \pc(\gh)} \theta(p) = 0$.

Consider the Bernoulli action of $\gp$ on $(X, \mu) :=
\big([0, 1]^\gp, \lambda^\gp\big)$, where
$\lambda$ is Lebesgue measure on $[0, 1]$.
The \dfn{Cayley graphing} of $X$ is the graph $(X, E_S)$, where $E_S :=
\big\{(x, s x) \st x \in X, s \in S \cup S^{-1}\big\}$.
For $x \in X$ and $p \in [0, 1]$, let $V_p(x) := \big\{\gpe^{-1} x \st \gpe
\in \gp, x(\gpe) \le p\big\}$.
Let $Y_p \subseteq X$ denote the set of points $x$ that belong to an infinite
cluster in the graph induced by $E_S$ on $V_p(x)$.
The orbit equivalence relation induces an equivalence relation $\rr_p$ on
$Y_p$.
Then $\mu(Y_p) = \theta(p)$ and for $p > \pc(\gh)$,
the cost of $\big(Y_p, \rr_p, (\mu \upharpoonright
Y_p)/\mu(Y_p)\big)$ is at most $|S|$. Hence by Proposition II.6 of \ref
b.Gaboriau:cost/, the cost of the Bernoulli action of $\gp$ on $(X, \mu)$
is at most $1 + \theta(p)\big(|S| - 1\big)$ for $p > \pc(\gh)$, whence the
cost equals 1.
Since this is the maximum cost of any action of $\gp$, while the minimum
cost is 1 for any infinite group, it follows that all costs are 1 and thus
that $\gp$ has fixed price.

We remark that one can replace the use of $\lim_{p \downarrow \pc(\gh)}
\theta(p) = 0$ by the fact that $\rr_{\pc}$ is hyperfinite, as pointed out
to us by D.~Gaboriau.
Indeed, the clusters induced by $E_S$ on
$V_p(x)$ for $p < \pc(\gh)$ are all finite by
definition for $\mu$-a.e.\ $x$, whence their increasing union is
hyperfinite and therefore has cost at most 1.
It follows that the cost of the Bernoulli action of $\gp$ on $(X, \mu)$
is at most $1 + \big(p - \pc(\gh)\big) |S|$ for $p > \pc(\gh)$, whence the
cost equals 1.

\medbreak
\noindent {\bf Acknowledgement.}\enspace 
I am grateful to Mikl\'os Ab\'ert for conversations.

\def\noop#1{\relax}
\input \jobname.bbl

\filbreak
\begingroup
\eightpoint\sc
\parindent=0pt\baselineskip=10pt

Department of Mathematics,
831 E. 3rd St.,
Indiana University,
Bloomington, IN 47405-7106
\emailwww{rdlyons@indiana.edu}
{http://mypage.iu.edu/\string~rdlyons/}

\endgroup


\def\cprime{$'$} \def\cprime{$'$} \def\cprime{$'$}
\def\temp{\let\linkit=\linkyear \apaliketrue}
\temp
\ifcitationgeneration\immediate\write\labelfile{\sanitize\temp}\fi
\def\startreferences{
 \vskip0pt plus.3\vsize \penalty -150 \vskip0pt
 plus-.3\vsize \bigskip\bigskip \vskip \parskip
 \begingroup\baselineskip=12pt\frenchspacing
 \bibliographytitle
 \vskip12pt\parindent=0pt
 \def\and{{\rm and}}
 \def\em{\it}
 \def\newblock{\hskip .11em plus.33em minus.07em}
 \def\bibauthor##1{{\sc ##1}}
 \def\bibitem[##1]##2
 {\htmlanchor{##2}{}\RefLabel{##2}[##1]\hangindent=.8cm\hangafter=1}
 }
\def\endreferences{\bigskip\bigskip\endgroup}
\ifundefined{bibstylemodification}\relax\else\bibstylemodification\fi
\startreferences

\bibitem[Ab\'ert and Weiss (2011)]{AW:Bernoulli}
\bibauthor{Ab\'ert, M. \and{} Weiss, B.} (2011).
\newblock Bernoulli actions are weakly contained in any free action.
\newblock Preprint, \arXiv{1103.1063}.

\bibitem[Benjamini, Lyons, Peres, and Schramm (1999)]{MR99m:60149}
\bibauthor{Benjamini, I., Lyons, R., Peres, Y., \and{} Schramm, O.} (1999).
\newblock Group-invariant percolation on graphs.
\newblock {\em Geom. Funct. Anal.} {\bf 9}, 29--66.

\bibitem[Benjamini and Schramm (1996)]{MR97j:60179}
\bibauthor{Benjamini, I. \and{} Schramm, O.} (1996).
\newblock Percolation beyond $\bold {Z}\sp d$, many questions and a few
  answers.
\newblock {\em Electron. Comm. Probab.} {\bf 1}, no.\ 8, 71--82 (electronic).

\bibitem[van~den Berg and Keane (1984)]{MR85g:60100}
\bibauthor{van~den Berg, J. \and{} Keane, M.} (1984).
\newblock On the continuity of the percolation probability function.
\newblock In Beals, R., Beck, A., Bellow, A., \and{} Hajian, A., editors, {\em
  Conference in Modern Analysis and Probability (New Haven, Conn., 1982)},
  pages 61--65. Amer. Math. Soc., Providence, RI.

\bibitem[Burton and Keane (1989)]{MR90g:60090}
\bibauthor{Burton, R.M. \and{} Keane, M.} (1989).
\newblock Density and uniqueness in percolation.
\newblock {\em Comm. Math. Phys.} {\bf 121}, 501--505.

\bibitem[Gaboriau (2000)]{MR1728876}
\bibauthor{Gaboriau, D.} (2000).
\newblock Co\^ ut des relations d'\'equivalence et des groupes.
\newblock {\em Invent. Math.} {\bf 139}, 41--98.

\bibitem[Gaboriau (2010)]{MR2761803}
\bibauthor{Gaboriau, D.} (2010).
\newblock What is {$\ldots$} cost?
\newblock {\em Notices Amer. Math. Soc.} {\bf 57}, 1295--1296.

\bibitem[Gandolfi, Keane, and Newman (1992)]{MR93f:60149}
\bibauthor{Gandolfi, A., Keane, M.S., \and{} Newman, C.M.} (1992).
\newblock Uniqueness of the infinite component in a random graph with
  applications to percolation and spin glasses.
\newblock {\em Probab. Theory Related Fields} {\bf 92}, 511--527.

\bibitem[Lyons (2000)]{MR2001c:82028}
\bibauthor{Lyons, R.} (2000).
\newblock Phase transitions on nonamenable graphs.
\newblock {\em J. Math. Phys.} {\bf 41}, 1099--1126.
\newblock Probabilistic techniques in equilibrium and nonequilibrium
  statistical physics.

\bibitem[Newman and Schulman (1981)]{MR83e:82038}
\bibauthor{Newman, C.M. \and{} Schulman, L.S.} (1981).
\newblock Infinite clusters in percolation models.
\newblock {\em J. Statist. Phys.} {\bf 26}, 613--628.

\bibitem[Pak and Smirnova-Nagnibeda (2000)]{MR1756965}
\bibauthor{Pak, I. \and{} Smirnova-Nagnibeda, T.} (2000).
\newblock On non-uniqueness of percolation on nonamenable {C}ayley graphs.
\newblock {\em C. R. Acad. Sci. Paris S\'er. I Math.} {\bf 330}, 495--500.

\endreferences
\bye